\documentclass[12pt]{article}
\usepackage{amssymb}
\usepackage{verbatim}
\usepackage{array}
\usepackage{latexsym}
\usepackage{enumerate}
\usepackage{amsmath}
\usepackage{amsfonts}
\usepackage{amsthm}
\usepackage{color}
\usepackage[english]{babel}

\title{{On the structure of $3$-nets \\ embedded  in a projective plane.}}
\date{}

\author{Aart~Blokhuis, G\'abor~Korchm\'aros, Francesco~Mazzocca}

\newtheorem{theorem}{Theorem}[section]

\newtheorem{lemma}[theorem]{Lemma}

{\theoremstyle{definition}
\newtheorem*{definition*}{Definition}

\newtheorem{rem}[theorem]{Remark}
\newtheorem*{proposition*}{Proposition}
\newtheorem*{corollary*}{Corollary}
\newtheorem*{lemma*}{Lemma}
\def\cC{\mathcal C}

\begin{document}
\maketitle
\begin{abstract}
We investigate finite $3$-nets embedded in a projective plane over a (finite or infinite) field of any characteristic $p$.
Such an embedding is {\em{regular}} when each of the three classes of the $3$-net comprises concurrent lines, and {\em{irregular}} otherwise. It is {\em{completely irregular}} when no class of the $3$-net consists of concurrent lines.
We are interested in embeddings of $3$-nets  which are irregular but the lines of one class are concurrent. For an irregular embedding of a $3$-net of order $n\geq 5$ we prove that,  if all lines from  two classes are tangent to the same irreducible conic,  then all lines from the third class are concurrent. We also prove the converse provided that the order $n$ of the $3$-net is smaller than $p$. In the complex plane, apart from a sporadic example of order $n=5$ due to Stipins \cite{sj2004}, each known irregularly embedded $3$-net has the property that all its lines are tangent to a plane cubic curve. Actually, the procedure of constructing irregular $3$-nets with this property works over any field. In positive characteristic, we present some more examples for $n\ge 5$ and  give a complete classification for $n=4$.
\end{abstract}
\section{The problem and the results}
A {\em {$3$-net of order $n$}} is a point-line incidence structure consisting of $n^2$ points together with three classes of lines each consisting of $n$ lines such that
\begin{itemize}
\item[\rm(i)] any two lines from different classes are incident;
\item[\rm(ii)] no two lines from the same class are incident;
\item[\rm(iii)] any point is incident with exactly one line from each class.
\end{itemize}

The notion of $3$-net comes from classical Differential geometry via the combinatorial abstraction of the notion of a $3$-web. There is a long history on $3$-nets in Combinatorics related to affine planes, latin squares, loops and strictly transitive permutation sets.  In recent years, finite $3$-nets embedded in the complex plane have been investigated in connection with the cohomology of local systems on the complements of complex line arrangements, see \cite{ys2004,fy2007,ys2007}. The methods exploit modern Topology and Algebraic geometry.

In the present paper, mostly combinatorial methods are used to investigate finite $3$-nets embedded in a projective plane $PG(2,{\bf F})$ over a field ${\bf F}$ of any  characteristic $p$ with special attention on finite projective planes.

For this purpose, it is convenient to regard $3$-nets in the dual plane of $PG(2,{\bf F})$; equivalently consider a dual $3$-net in $PG(2,{\bf F})$. Here {\em{a dual $3$-net embedded in}} $PG(2,{\bf F})$ is a triple $\{A,B,C\}$ with $A,B,C$ pairwise disjoint point-sets of size $n$, called {\em{components}}, such that every line meeting two distinct components meets each component in precisely one point. A dual $3$-net is {\em{regularly}} embedded if each component is contained in a line of $PG(2,{\bf F})$, otherwise the embedding is {\em{irregular}}. It is {\em{completely irregular}} if no component is contained in a line. Many examples of completely irregular embeddings are given by the dual $3$-nets whose points lie on the same irreducible plane cubic curve.
Our main results are stated in Theorems \ref{thm1} and \ref{AUBconic}. They provide a characterization of dual $3$-nets with a non-completely irregular embedding in $PG(2,{\bf F})$.


We also show that if $n=4$ the set $A\cup B\cup C$ is contained in a cubic, this is no longer true in general when $n\ge 5$.

\section{Dual $3$-nets of small order}

For $n=1$, any dual $3$-net of order $n$ is trivial as each component is just a point.
For $n=2$, for every dual $3$-net $\{A,B,C\}$  in $PG(2,{\bf F})$, $A\cup B\cup C$ is projectively equivalent to the \it Pasch configuration. \rm Under suitable choice of the coordinate system, the  Pasch configuration consists of the points
$$(0,0,1),\,(0,1,0),\,(1,0,0),\,(1,1,1),\,(1,0,1),\,(0,1,1),$$
lying on each of the plane cubic curves belonging to the linear system  $$aX^2Y+bXY^2+c(X^2Z-XZ^2)+d(Y^2Z-YZ^2)-(a+b)XYZ.$$

For $n=3$, all irregular dual $3$-nets are projectively equivalent to $\{A,B,C\}$ where
$$
\begin{array}{lll}
A &=& \{(1:0:0),\,(0:1:0),\,(0:0:1)\},\\
B &=& \{(a:1/b:1),\,(b:1/c:1),\,(c:1/a:1)\},\\
C &=& \{(a:1/c:1),\,(b:1/a:1),\,(c:1/b:1)\},
\end{array}
$$
with three pairwise distinct non-zero elements $a,b,c$
of ${\bf F}$. Since a dual $3$-net of order $3$ consists of nine points, they are contained in a cubic.
Furthermore, $B$ is on a line when
$$\frac{a}{c}+\frac{c}{b}+ \frac{b}{a}=3,$$
and similarly for $C$ when
$$\frac{a}{b}+\frac{b}{c}+ \frac{c}{a}=3.$$
From this, if $B$ is on a line, then $A\cup C$ is on an irreducible conic with equation of type
$\alpha XY +\beta YZ +\gamma ZX=0$ with $\alpha,\beta,\gamma\in {\bf F}$.
The converse is also true, that is, if $A\cup C$ is on an irreducible conic then $B$ is on a line.
The last two results are special cases of Theorems \ref{thm1} and \ref{AUBconic}.

The investigation of the case $n=4$ requires much effort and it is carried out in Section \ref{sec4}.

\section{Examples of dual $3$-nets from cubic curves}
\label{excubic}
We show a general construction for dual $3$-nets from a plane cubic curve ${\cal F}$ of
$PG(2,{\bf F})$. The procedure is a straightforward extension of that used in \cite{ys2004} and it depends on the well known abelian group $(G,+)$ defined on the points of a  cubic curve.

\begin{theorem}{\rm{\cite[Theorem 6.104, Corollary 6.105]{hkt}}}
\label{groupcubic} A non-singular plane cubic curve ${\cal F}$ can be equipped with an operation on its points to obtain an abelian
group $(G,+)$. If  $0$ is the identity of $(G,+),$ then three distinct
points $P,Q,R\in {\cal F}$ are collinear if and only if $P+Q+R=0',$ $0'$ being the third of the three points $($counted with multiplicity$)$ that the tangent line in $0$ to $\cal F$ has in common with $\cal F.$
\end{theorem}

For sake of completeness, we recall that this group is defined as follows. Choose a point $O\in \cal F$ as the identity of the group. For any two points $P,Q\in\cal F ,$ the sum of $P$ and $Q$ is defined by $P+Q=R, $ $R$ being the third of the three points (counted with multiplicity) that the line through $P,Q$ has in common with $\cal F.$
In the same way a group can be defined on the set of non singular points when the cubic splits into an irreducible conic and a line [see the final remark in Section \ref{oneset}].

Let $G$ be the abelian group associated to a non-singular cubic curve ${\cal F}$  in Theorem \ref{groupcubic}. Take a  finite subgroup $H$ of $G$ whose index is greater than two, with $0'\in H,$ and choose three
pairwise distinct cosets of $H$ in $G$, say
$$ A=a+H,\, \, B= b+H,\, \, C=c+H,  $$ with $a,b,c\in G$ and collinear, i.e. $$  a+b+c=0' \, . $$
Then $A\cup B \cup C$ is a dual $3$-net whose order is equal to the size of $H$. If $0$ is an inflection point of  ${\cal F}$, then $0=0'$ and the collinearity condition for the points $a,b,c$ becomes  $a+b+c=0$. A $3$-net arising in this way from a cubic curve is said to be {\em of subgroup type} if $0$ is an inflection point.

The above construction, in case $0$ is an inflection point, is currently the main source of examples of $3$-nets. Even in finite planes, it provides a lot of examples using the following result on the number of points of a non-singular cubic curve defined over a finite field.
\begin{theorem}{\rm\cite[Waterhouse's theorem, Theorem 9.92]{hkt}}
Let ${{\bf{F}}_q}$ be a finite field of order $q$. For every positive integer
$m$ such that $|m|\leq 2\sqrt{q}$ and $m\not\equiv 0 \pmod p$, there exists a non-singular cubic curve ${\cal F}$ defined over ${{\bf{F}}_q}$ with precisely $N=q+1-m$ points in $PG(2,{\bf{F}}_q)$.
\end{theorem}
Here ${\cal F}$ may be assumed to be a non-singular plane cubic curve with an inflexion in $PG(2,{\bf{F}}_q)$,  see \cite[Chapter 9.10]{hkt}).

The known results on finite $3$-nets embedded in the complex plane mentioned in Introduction and the above discussion might suggest the conjecture that
the completely irregular finite dual $3$-nets, apart from a few mostly sporadic examples, should be contained in a cubic curve.  We have already quoted Stipins' sporadic example of order $n=5$ in the complex plane.
In $PG(2,{\bf F}_8)$ there is a nice example of order $n=5$ due to J. Bierbrauer \cite{bb}. In this example, the point-sets $A\cup B$, $B\cup C$ and $C\cup A$ are hyperovals, and $A\cup B\cup C$ is a 3-arc of maximal size not contained in a cubic curve.



\section{Examples with the three sets not contained in the same cubic curve}
\label{sec5}
An irregular dual $3$-net with a component on a line is contained in a cubic curve if and only if the other components lie on an irreducible conic.
This raises the question of determining the irregular dual $3$-nets such that a component is on a line but there is no irreducible conic containing the other two components.

We show that such irregular dual $3$-nets exist, indeed. They can even be constructed in a finite plane $PG(2,{\bf{F}}_q)$ provided that ${\bf F}_q$ has a subfield ${\bf F}_r$ of order $r>3$ satisfying the arithmetic condition
\begin{equation}
\label{eqq1} q>r^2.
\end{equation}
To show how such a construction can be performed, consider the three-dimensional affine space
$AG(3,{\bf{F}}_r)$, and embed it canonically in the three-dimensional affine space $AG(3,{\bf{F}}_q)$.
Let $n=r^2$ and choose three pairwise distinct parallel planes $\alpha,\beta$ and $\gamma$ in $AG(3,{\bf{F}}_r)$.
Obviously,
\begin{itemize}
\item[(*)]
any line meeting two of these planes must meet the third one too, and in this case the line meets each these planes once.
\end{itemize}
Consider the unique plane $\gamma'$  of $AG(3,{\bf{F}}_q)$ containing $\gamma$, and
take a point $P$ so that no line containing more than one point of $\gamma$ contains $P$. This can be done since $\gamma$
can be viewed as an affine subplane of $\gamma'$ and its order $r$ satisfies (\ref{eqq1}).
Now, in the three-dimensional projective space $PG(3,{\bf{F}}_q)$ arising from $AG(3,{\bf{F}}_q)$, project the points on $\alpha,\beta$ and $\gamma$ from $P$ on a plane $\pi$ disjoint from $\alpha,\beta,\gamma$ and not passing through $P$.  The projection takes $\alpha,\beta$ and $\gamma$ to three point-sets $A,B$ and $C$ in $\pi$ each of size $n=r^2$. Obviously, $C$ is contained in the common line of $\pi$ and $\alpha'$. We prove the following properties:
\begin{itemize}
\item[{\rm(I)}] $A$ and $B$ have no common point and they are also disjoint from $C$.
\item[{\rm(II)}] Any line joining a point of $A$ with a point of $B$ meets $C$, and it meets each of $A,B$ and $C$ only once.
\end{itemize}
Assume on the contrary that $A\cap B$ contains a point $Q$. Then the points $\bar{A}=PQ\cap \alpha$ and $\bar{B}=PQ\cap\beta$ are two distinct points in $AG(3,{\bf F}_r)$. This implies that $P$ is on a line $\ell$ of $AG(3,{\bf F}_r)$, and hence $P$ is the common point of $\ell$ with the plane $\gamma'$. But this is impossible, since $P$ is not in $\gamma$.

 To show (II), take any two points $A'\in A$ and $B'\in B$. Then $A'\neq B'$ by (I). Let $\bar{A}=PA\cap \alpha$ and $\bar{B}=PB\cap\beta$. Then the line $\bar{\ell}$ joining $\bar{A}$ with $\bar{B}$ is in $AG(3,{\bf{F}}_r)$ and hence it meets
$\gamma$ in a point $\bar{C}$. Now, the projection from $P$ takes $\bar{C}$ to a point in $C$ which lies on $\ell$. Finally, assume on the contrary that a line joining two points $A_1,A_2\in A$ contains either a point $B'\in B$, or a point $C'\in C$. To deal with the former case, let $\bar{A}_1=PA_1\cap \alpha,\bar{A}_2=PA_2Q\cap \alpha$ and $\bar{B}=PB'\cap \beta$. Then the plane $\delta$ determined by the triangle $\bar{A}_1\bar{A}_2\bar{B}$ passes through $P$. Therefore, $P$ is on the common line $\ell$ of two planes $\delta$ and $\gamma$ of $AG(3,{\bf{F}}_r)$. But this contradicts the choice of $P$ as $\ell$ contains at least two points from $\gamma$. In the latter case, define $\bar{A}_1$ and $\bar{A}_2$ as before, and let $\bar{C}=PC'\cap \gamma$. Arguing as before we get a contradiction as the common line $\ell$ of $\gamma$ and the plane $\delta$ determined by the triangle $\bar{A}_1\bar{A}_2\bar{B}$ passes through $P$. This completes the proof of (II).

Therefore, $\{A,B,C\}$ is a dual $3$-net embedded in the plane $\pi\cong PG(2,{\bf F}_q)$ such that $C$ is on a line. Since $A$ and $B$ are affine planes of order $r$, they contain $r^2$ subsets of $r$ collinear points but they do not contain $r+1$ points from a line. Since $r>3$, neither $A$ nor $B$ are contained in a plane cubic curve, and the same holds for $A\cup B$.
Our final remark is that the above construction may be modified in such a way to produce a dual $3$-net with a complete irregular embedding in $\pi$ so that no component is contained in a
(possibly reducible) plane cubic curve.

\section{One component is contained in a line}
\label{oneset}
As we have seen in Section \ref{sec5}, there exist dual $3$-nets irregularly embedded in some projective plane $PG(2,{\bf F})$ such that a component is on a line but there is no irreducible conic containing the other two components. In this section, we prove however that this does not happen in any $PG(2,{\bf F})$, in particular for fields ${\bf F}$ of zero characteristic.

\begin{theorem}
\label{thm1}
Let $\{A, B, C\}$ be a dual $3$-net of order $n$ in $PG(2,{\bf F})$. When ${\bf F}$ has positive characteristic $p$ suppose that $n\leq p$.
If $C$ is on a line, then $A\cup B$ is on a conic.
\end{theorem}

For the line containing $C$ we take $Z=0$, the line at infinity. Then $A$ and $B$ will be in $AG(2,{\bf F})$.

Two points $a=(a_1,a_2)$ and $b=(b_1,b_2)$ in $AG(2,{\bf F})$ determine
the direction
\[
m={b_2-a_2\over b_1-a_1}
\]
This we can rewrite as $mb_1-b_2=ma_1-a_2$.

Consider a transversal: $Y=mX+f$, containing $a=(a_1,a_2)\in A$ and $b=(b_1,b_2)\in B$ (and $(1:m:0)=:(m)$ in $C$)
then $ma_1-a_2=mb_1-b_2$. Consider two R\'edei polynomials (of degree $n$)
\[
A(T,X)=\prod_{a\in A} (T+(Xa_1-a_2));\quad B(T,X)=\prod_{b\in B} (T+(Xb_1-b_2)).
\]
If $X=m$ and $(m)\in C$ then the polynomials $A(T,m)$ and $B(T,m)$ are the same.
This means that $A(T,X)$ and $B(T,X)$ are
almost identically the same, because all coefficients of their difference are divisible by
the polynomial $\prod_{(m)\in C} (X-m)$, of degree $n$. The only coefficient of high enough degree
(in $X$) to be nonzero, is the 'constant term' that is, the term not containing $T$, of $A-B$.
So $\prod (Xa_1-a_2) - \prod (Xb_1-b_2)$ is divisible
by (and in fact a scalar multiple of) $\prod_{(m)\in C} (X-m)$.
Here we tacitly assumed that the vertical (infinite) direction does not occur, this of
course we may do without loss of generality.\\
We have that the elementary
symmetric polynomials $\sigma_k(A)$ (by this we mean the coefficient of $T^{n-k}$ in $A(T,X)$)
and $\sigma_k(B)$ are equal for $k=0,1,\dots,n-1$.
Hence also the power sum polynomials $\pi_k(A):=\sum_{a\in A} (Xa_1-a_2)^k$ and $\pi_k(B)$ are equal for these $k$.
This means that for all $i,j$ with $i+j<n$: $\sum_{a\in A} a_1^ia_2^j=\sum_{b\in B} b_1^ib_2^j$ by comparing
the coefficient of $X^i$ in $\pi_{i+j}$ (here we use the condition that $n\le p$, if the characteristic is positive). This now means that if we have {\em any}
polynomial $f(u,v)$ in two variables of total degree at most $n-1$, then
$\sum_{a\in A} f(a)=\sum_{b\in B} f(b)$. This then implies that there is no polynomial of
total degree at most $n-1$ that vanishes on all but one of the points of $A\cup B$, for then one of the sides of this equation vanishes, and the other doesn't. Now this
indeed is the case if both $A$ and $B$ consist of $n$ collinear points, or, more generally if $A\cup B$ is contained in a conic (by Bezout's theorem), which we are now going to
show. Any 5 points (let's say 3 from $A$ and 2 from $B$) determine a conic,
and if this conic does not contain all points of $A\cup B$, then we can cover the remaining
points except one with lines (this is because any line containing a point from $A$ and
a point from $B$ contains no further points of $A$ or $B$).
In this way we find a polynomial $f$ of total
degree at most $n-1$ vanishing on all but one of the points of $A\cup B$ and thus get a contradiction. So $A\cup B$ lies on a conic.\\
The special situation that $A$ and $B$ are contained in a line is of course well known.
To make the paper self contained we consider it here in more detail. If the lines
are not parallel, we may take $A$ on the $X$-axis, $B$ on the $Y$-axis, and then the
points $(a,0)$ and $(0,b)$ will determine the direction $m=-b/a$. It follows that $A$ and $B$ are both cosets of a subgroup $G$ of the multiplicative group ${\bf F}^*$.\\
If the lines are parallel, say $X=0$ and $X=1$, then $(0,a)$ and $(1,b)$ determine the
direction $m=b-a$ and we get that $A$ and $B$ are cosets of an additive subgroup G
of ${\bf F}$. If the characteristic is zero this situation does not occur, since there
are no finite subgroups. If the characteristic is $p$ and $n\le p$ then $G$ is necessarily generated by a single element. Of course for $n>p$ we get more examples.\\
Next we consider the case that $A\cup B$ is contained in an irreducible conic ${\cal C}$.
Since the line at infinity is special we get three cases. First we investigate the case
that ${\cal C}$ is a parabola, which we take to be $Y=X^2$.
Let $A$ consists of pairs $(a,a^2)$, $B$ of pairs $(b,b^2)$,
so the directions will be $(b^2-a^2)/(b-a)=b+a$ so $A$ and $B$ correspond to cosets
of some subgroup $G$ of the additive group of the field. Next we look at the case that
$\cal C$ is the hyperbola $XY=1$, $A$ and $B$ having points
$(a,1/a)$ and $(b,1/b)$ determining the direction $(1/b-1/a)/(b-a)=-1/ab$.
Now $A$ and $B$ are cosets of $G$, a subgroup of ${\bf F}^*$.\\
Finally consider a conic contained in the affine plane. Let ${\bf F}=GF(q)$ and identify $AG(2,q)$ with $GF(q^2)$.
Now let ${\cal C}$ be the 'circle' $U^{q+1}=1$. The points on this circle form the cyclic group $GF(q^2)^{*(q-1)}$ of order
$q+1$. Two points $a,b\in {\cal C}=GF(q^2)^{*(q-1)}$  determine direction
\[
(a-b)^{p-1}={a^p-b^p\over a-b}={1/a-1/b\over a-b}={-1\over ab}.
\]
So if the set $\{ab\,|\,a\in A,b\in B\}$ has size $n$, then after a multiplicative shift we
may assume that 1 is in both sets $A$ and $B$ and both must be the same subgroup (of the multiplicative group of
$q+1$-st roots of unity).
\begin{rem}
\label{rmk}
A more unified approach to the above is the geometric description of the group $G$, in the case where a plane cubic splits into an irreducible conic and a line.
So let ${\cal C}$ be an irreducible conic and $\ell$ a line. Let $O$ be a point on $\cC$ but on $\ell$. Then an abelian group on the points on $\cC$ not lying on $\ell$ arises from the operation
$a+b=c$ defined geometrically as follows.  Let $m$ be the chord of $\cC$ through $a$ and $b$ when $a\neq b$ and the tangent to $\cC$ at $a$ when $a=b$. Let $P$ be the point cut out on $\ell$ by $m$.
If $OP$ is chord of $\cC$, then $c$ is defined to be point cut out on $\cC$ by $OP$. If $PO$ is the tangent to $\cC$ at $O$, then $a+b=O$. Now, if the point $O$ is taken from $A$, then
$A + B = B$, whence  $A + A + B = B$. Therefore, $A$ is a subgroup, and $B$ is a coset of $A$, and $\{A,B,C\}$ is of 'subgroup type'.
\end{rem}

\section{Two components are on the same irreducible conic}
Our goal is to prove the converse of Theorem \ref{thm1}. For $n=4$ this follows from Theorem \ref{thm2}.
\begin{theorem}
\label{AUBconic}
Let $\{A, B, C\}$ be a dual $3$-net of order $n\geq 5$ in $PG(2,{\bf F})$ such that
$A\cup B$ is contained in a conic. Then $C$ is contained in a line.
\end{theorem}
Let $\cC$ be the conic containing both $A$ and $B$. Without loss of generality, we may assume that $\cC$ is irreducible. In fact, if $\cC$ splits into two lines, then
one of them contains $A$ and the other does $B$. Therefore, the assertion follows from Theorem \ref{thm1}.

{}From now on, $\cC$ is assumed to be irreducible. The key idea of the proof is to project $A$ into $B$ from a point $Q\in C$ and investigate the permutation group $\Phi$ on $A\cup B$ generated by all such projections. The hypothesis that $A\cup B$ is contained in $\cC$ ensures that $\Phi$ is isomorphic to a finite subgroup of $PGL(2,{\bf F})$. Using the classification of finite subgroups  $PGL(2,{\bf F})$ we will be able to determine the possibilities for $\Phi$. Actually, only one possibility can occur, namely when $\Phi$ is a dihedral group. On the other hand, for a dihedral group $\Psi$, we will prove that the points in $C$ must be on a line which is either a secant or an external line to $\cC$.

For the proofs, some preliminary results on linear collineations of an irreducible conic are collected.

The full linear collineation group $G$ of $PG(2,{\bf F})$ which preserves $\cC$ is isomorphic to $PGL(2,{\bf F})$, and $G$ acts on $\cC$ as $PGL(2,{\bf F})$ on the projective line $PG(1,{\bf F})$. If $$R=:\,\,x\mapsto\frac{ay+b}{cy+d}\quad a,b,c,d\in {\bf F}\,\,with\,\,ad-bc\neq 0$$ is an element from $PGL(2,{\bf F})$
then the corresponding element of $G$ is the linear collineation associated with the non-singular matrix
$$R'=
\left(
  \begin{array}{ccc}
    a^2 & 2ab & b^2 \\
    ac & ad+bc & bd \\
    c^2 & 2cd & d^2 \\
  \end{array}
\right),
$$
see \cite[Theorem 2.37]{hp1973}. Here, $\cC$ consists of all points $P_y(y^2,y,1)$ together with $P_\infty=(1,0,0)$. Furthermore, $R'$ takes $P_\infty$ to $(a^2,ac,c^2)$ while it takes $P_y$ to either $P_\infty$ or $$(\left(\frac{ay+b}{cy+d}\right)^2,\frac{ay+b}{cy+d},1)$$ according as $cy+d$ vanishes or does not.

Also, $R'$ has $0,1$ or $2$ fixed points and when $R$ has no fixed point then it has two fixed points over the quadratic extension ${\bf F'}$ of ${\bf F}$. This follows from the corresponding result on the action of $R$ on $PG(1,{\bf F})$.

For a point $Q$ outside $\cC$ (and also different from the nucleus of $\cC$ when $p=2$), the projection of $\cC$ onto itself from $Q$ is the restriction on $\cC$ of unique  involutory perspectivity $\varphi_Q\in G$. The center of $\varphi_Q$ is $Q$. So, two distinct points $C_1,C_2\in \cC$ correspond under the action of $\varphi_Q$ if and only if  the line through $C_1$ and $C_2$ also contains $Q$. 

From now on $\Phi$ denotes the subgroup of $G$ generated by all involutory perspectivities $\varphi_Q$ with $Q$ ranging over $C$. Obviously, every $\varphi_Q$ interchanges $A$ and $B$, so that the product of any two such  involutory perspectivities preserves both $A$ and $B$. This suggests to consider the subgroup $\Psi$ of $\Phi$ generated by all products $\varphi_Q\varphi_S$ with $Q,S\in A$. Since $A$ and $B$ are disjoint, $\Psi$ has index two in $\Phi$. Some more properties of $\Phi$ are given in the following result.
\begin{lemma}
\label{lem1}
$\Psi$ acts transitively on both $A$ and $B$. If $\Psi$ is also regular on $A$ then $\Psi$ is an abelian group of order $n$.
\end{lemma}
To show the first assertion, let $A_1,A_2\in A$ be two distinct points. Take any point $V\in B$, and consider the points $Q,S\in C$ cut out on $C$ by the lines $VA_1$ and $VA_2$, respectively. Then $\varphi_Q\varphi_S$ takes $A_1$ to $A_2$, and the assertion follows. Now suppose that $\Psi$ is regular on $A$. Then $\Psi$ has order $n$. Hence, for a given $Q\in C$, every element in $\Psi$ is of the form $\varphi_Q\varphi_S$. Furthermore,
$\Phi$ has order $2n$. Therefore, the coset $\Phi\setminus \Psi$ consists of $n$ elements. On the other hand, each $\varphi_Q$ with $Q\in C$ is in that coset. Since the number of such $\varphi_Q$ is $n$, it turns out that the coset entirely consists of involutions. Therefore, the product $\varphi_Q \varphi_R\varphi_S$ with $Q,R,S\in A$ is also
an involutory perspectivity $\varphi_T$ for some $T\in C$. In particular, $(\varphi_S\varphi_Q\varphi_T)^2=1$. From this, $(\varphi_Q\varphi_S)(\varphi_Q\varphi_T)=(\varphi_Q\varphi_T)(\varphi_Q\varphi_S)$ for any two $S,T\in A$. Therefore $\Psi$ is abelian, and the proof of Lemma \ref{lem1} is complete.

We remark that $\Psi$ is not regular if and only if the stabilizer of a point $P\in A$ of $\Psi$ is non-trivial. If this is the case, then $A\cup B$, viewed as an orbit of $\Phi$, is a short orbit, that is, its size is smaller than the order of $\Phi$. 

Now, the classification theorem of all finite subgroups of $PGL(2,{\bf F})$ is stated. For a proof, see \cite{mv1983}.
\begin{theorem}
\label{vm1} Let $U$ be a non-trivial finite
group preserving an irreducible conic $\cC$ in $PG(2,{\bf F})$.
Let $s$ be the number of short orbits of $U$ on the set of all
points on the irreducible conic $\cC'$ which is the extension of $\cC$ in $PG(2,{\bf F'})$ where ${\bf F'}$ is a degree two algebraic extension of ${\bf F}$.
Let $\ell_1,\ldots,\ell_s$ be the lengths of the short orbits of $G$ on $\cC'$.
Then $U$ is a group of one of the following types$:$
\begin{enumerate}
\item[\rm(i)] the cyclic group $C_m$ of order $m,$ with $p\nmid m,
\,s=2,\, \ell_1=\ell_2=1;$

\item[\rm(ii)] an elementary abelian $p$-group$,$ with $s=1,\,
\ell_1=1;$

\item[\rm(iii)] the dihedral group $D_m$ of order $2m,$ with
$p\nmid m,\, s=2,\, \ell_1=2,\, \ell_2=m,$ or $p\neq 2,\, s=3,\,
\ell_1=\ell_2=m,\,\ell_3=2;$

\item[\rm(iv)] the alternating group ${\bf A} _4,$ with $p\neq 2,3,$ and
$s=3,\, \ell_1=6,\, \ell_2=\ell_3=4;$

\item[\rm(v)] the symmetric group ${\bf S}_4,$ with $p\neq 2,3,\, s=3,\,
\ell_1=12,\, \ell_2=8,\, \ell_3=6;$

\item[\rm(vi)] the alternating group ${\bf A}_5,$ with $p=3,\, s=2,\,
\ell_1=10,\, \ell_2=12,$ or\linebreak $p\neq 2,3,5,\, s=3,\, \ell_1=30,
\,\ell_2=20,\, \ell_3=12;$

\item[\rm(vii)] the semidirect product of an elementary abelian
$p$-group of order $q$ with a cyclic group of order $m,$ with
$m\mid(q-1),\, s=2,\, \ell_1=1,\, \ell_2=q;$

\item[\rm(viii)] $PSL(2,q),$ with $p\neq 2,\, q=p^m,\, s=2,\,
\ell_1=q(q-1),\, \ell_2=q+1;$

\item[\rm(ix)] $PGL(2,q),$ with $q=p^m,\, s=2,\,
\ell_1=q(q-1),\, \ell_2=q+1$.
\end{enumerate}
\end{theorem}

We are in a position to prove Theorem \ref{AUBconic}. To apply Theorem \ref{vm1}, $\Phi$ has to be considered as a subgroup of $PGL(2,{\bf F'})$.

Consider at first the possibility that $\Psi$ is not regular on $A$. Then $A$ and $B$ are short orbits of $\Phi$.
Since $A$ and $B$ have the same size $n$ with $n\geq 5$, Theorem \ref{vm1} implies that $\Psi$ is a dihedral group $D_n$ of order $2n$ with $p\nmid n$. Then $\Psi$ has a short orbit of size $2$ in $\cC$ or in its extension $\cC'$. Let $P_1,P_2$ be the two distinct points which form that short orbit. Theorem \ref{vm1} also implies that the unique cyclic subgroup of $\Phi$ fixes both $P_1$ and $P_2$. Therefore, every $\varphi_Q$ interchanges $P_1$ and $P_2$.
Consider the line $\ell$ through $P_1$ and $P_2$. Even though $P_1$ and $P_2$ may only defined over in $PG(2,{\bf F'})$,  $\ell$ is line in $PG(2,{\bf F})$. To show that $A$ is contained in $\ell$ take a point $Q\in A$.
Then $\varphi_Q$, viewed as an  involutory perspectivity in $PG(2,{\bf F'}$ when $P_1,P_2$ are defined over ${\bf F'}$, interchanges $P_1$ and $P_2$. Therefore, $\ell$ contains the center of the involutory perspectivity $\varphi_Q$. Since $Q$ is the center of $\varphi_Q$, the assertion follows.

Assume that $\Psi$ is regular on $A$. From Lemma \ref{lem1} and Theorem \ref{vm1}, only two possibilities may occur for $\Phi$, namely (iii), and (vii) for $m=2$. In the former case, $\Phi$ is a dihedral group and the above argument involving the short orbit $\{P_1,P_2\}$ still works and it shows that $A$ is contained in the line $\ell$. In the latter case, $\Psi$ comprises elements of order $p$ together with the identity. This yields that $\Psi$ fixes a unique point $P\in\cC$. Hence, $\Psi$ preserves the tangent line $\ell$ to $\cC$ at $P$. Since $\Psi$ is a normal subgroup of $\Phi$, this implies that $\Phi$ itself preserves $\ell$. If $\varphi_Q$ fixes $\ell$ pointwise then $p=2$ and $\ell$ is the axis of $\varphi_Q$.  Therefore $Q\in \ell$ in this case, since the axis and the center of an involutory perspectivity are coincident  when $p=2$. Otherwise, $p\neq 2$ and $\ell$ is not the axis of $\varphi_Q$. Hence $\ell$ contains the center of $\varphi_Q$ which is the point $Q$. This completes the proof of Theorem \ref{AUBconic}.

\section{The case $n=4$}
\label{sec4}
 In the examples presented in Section \ref{sec5}, $n$ is always larger than the characteristic $p$.  More precisely
 $n\ge p^2$, so we could still think that the conjecture stated in Section \ref{excubic} might hold
 for smaller $n$ or in characteristic zero.
As a matter of fact, in the example by J. Bierbrauer, if $n>4$, the set
$A\cup B\cup C$ is not necessarily contained in a cubic. However, for $n=4$ we have the
following result.
\begin{theorem}
\label{thm2}
Let $\{A, B, C\}$ be a dual $3$-net of order 4 in $PG(2,{\bf F})$. Then $A\cup B\cup C$
is contained in a cubic.
\end{theorem}

In fact we can completely characterize the case $n=4$. If $A$, $B$ and $C$ are
contained in a line, then clearly their union is contained in the cubic that is the
product of these three lines, and we have one of the examples described in the previous
section. Next we consider the case that the set $A$ contains three non collinear points: $a_1=(1:0:0)$, $a_2=(0:1:0)$, $a_3(0:0:1)$ and a fourth point which is either
$(1:1:1)$ (case A, for arc) or $(1:1:0)$ (case non-A). Let $B=\{(a:b:1)$, $(c:d:1)$, $(e:f:1)$, $(g:h:1)\}$. Note that
$B$ cannot have a point on the line at infinity $Z=0$, because that line contains
two points of $A$ so we may normalize the third coordinate of each point to 1 and now
it follows that $a,b,c,d$, $e,f,g$ and $h$ are all different and nonzero, for otherwise
there would be a line containing at least two points of $B$ and a point of $A$.
For $C$ we have (up to reordering the points) only two possibilities: $C_1=\{(a:d:1)$, $(c:f:1)$, $(e:h:1)$, $(g:b:1)\}$ (case C, for cyclic), or $C_2=\{(a:d:1)$, $(c:b:1)$, $(e:h:1)$, $(g:f:1)\}$, (non-C). Why is this?
Because of the point $a_2\in A$ the first coordinates of the points in $C$ are the same as those in
$B$. Because of the point $a_1$ the second coordinates are also the same.  There are only 2
conjugacy classes of
permutations of 4 symbols without fixed points, a 4-cycle or the product of two 2-cycles.
Now the point $a_3\in A$ implies that the set of ratios $\{a/b,c/d,e/f,g/h\}$ is the same as $\{a/d,c/f,e/h,g/b\}$ in case C, and $\{a/d,c/b,e/h,g/f\}$ in case non-C. In the first case we have without loss of generality
$a/b=c/f$; $c/d=e/h$; $e/f=g/b$ and $g/h=a/d$, in other words
\[
af=bc, ~~~ ch=de, ~~~ eb=gf, ~~~ gd=ah.
\]
In the second case we have (again without loss of generality) $a/d=e/f$; $c/b=g/h$; $e/h=a/b$ and $g/f=c/d$:
\[
af=de, ~~~ ch=bg, ~~~ eb=ah, ~~~ gd=fc.
\]
Another way to see that there are essentially these two possibilities is to realize that
to a $3$-net there corresponds a latin square (by letting $A,B$ and $C$ play the
r\^ole of rows, columns and symbols), and there are only two essentially different latin squares of order 4.\\
We start to investigate the first case, where we will separate the cases A and non-A later.

To show that the points are on a cubic we consider the table
\[
\begin{array}{l|rrrrrrr}
P     & x^2y & x^2z & y^2x & y^2z & z^2x & z^2y & xyz\\ \hline
(ab1) & a^2b & a^2  & b^2a & b^2 &  a   &  b   & ab\\
(cd1) & c^2d & c^2  & d^2c & d^2 &  c   &  d   & cd\\
(ef1) & e^2f & e^2  & f^2e & f^2 &  e   &  f   & ef\\
(gh1) & g^2h & g^2  & h^2g & h^2 &  g   &  h   & gh\\
(ad1) & a^2d & a^2  & d^2a & d^2 &  a   &  d   & ad\\
(cf1) & c^2f & c^2  & f^2c & f^2 &  c   &  f   & cf\\
(eh1) & e^2h & e^2  & h^2e & h^2 &  e   &  h   & eh\\
(gb1) & g^2b & g^2  & b^2g & b^2 &  g   &  b   & gb
\end{array}
\]
and prove that the corresponding matrix has rank less than 7. After some row operations
(and using the fact that all letters are different nonzero field elements) we obtain:
\[
\left(
\begin{array}{rrrrrrr}
 0 & a^2  & -bad & -bd &  a   &  0   &  0\\
 0 & c^2  & -dcf & -df &  c   &  0   &  0\\
 0 & e^2  & -feh & -fh &  e   &  0   &  0\\
 0 & g^2  & -hgb & -hb &  g   &  0   &  0\\
 a^2  &   0  & ad+ab & d+b &  0   &  1   & a\\
 c^2  &   0  & cf+cd & f+d &  0   &  1   & c\\
 e^2  &   0  & eh+ef & h+f &  0   &  1   & e\\
 g^2  &   0  & gb+gh & b+h &  0   &  1   & g
\end{array}
\right)
\]
After some more manipulations we get
\[
\left(
\begin{array}{rrrrrrr}
 0 & a^2  & -bad & -bd &  a   &  0   &  0\\
 0 & c-a  & -df+bd & 0 &  0   &  0   &  0\\
 0 & e-a  & -fh+bd & 0 &  0   &  0   &  0\\
 0 & g-a  & -hb+bd &  0 &  0   &  0   &  0\\
 a^2  &   0  & ad+ab & d+b &  0   &  1   & a\\
 c^2  &   0  & cf+cd & f+d &  0   &  1   & c\\
 e^2  &   0  & eh+ef & h+f &  0   &  1   & e\\
 g^2  &   0  & gb+gh & b+h &  0   &  1   & g
\end{array}
\right)
\]
The second, third and fourth  row of this matrix
are scalar multiples of each other. So for a dependency
relation we may take $(x_1,x_2=\lambda hb,x_3=\lambda g,x_4,x_5=(bd/a)x_4,x_6,x_7)$,
where $(x_1,x_4,x_6,x_7)$ is a solution of the matrix
equation $A{\bf x}={\bf b}$ with $(A|{\bf b})$ given by:
\[
\left(
\begin{array}{rrrr|r}
 a^2  & d+b &  1   & a & -\lambda ag(d+b)\\
 c^2  & f+d &  1   & c & -\lambda cg(f+d)\\
 e^2  & h+f &  1   & e & -\lambda eg(h+f)\\
 g^2  & b+h &  1   & g & -\lambda g^2(b+h)
\end{array}
\right)
\]
The solution we find, with suitable $\lambda$, is:
\begin{eqnarray*}
x_1 &=& bf- dh,\\
x_2 &=& -bfd-hbf+hbd+hdf,\\
x_3 &=& -bc-be+ah+de,\\
x_4 &=& ae( b+f-d-h),\\
x_5 &=& bde(b+f-d-h),\\
x_6 &=& ae(-bf+dh),\\
x_7 &=& (bd-fh)(c-g).
\end{eqnarray*}
Finally we consider the two possibilities for the fourth point of $A$.
If we are in case A, then the fourth point is $(1:1:1)$ and in this case
we must have $\sum x_i=0$, where we have the additional equations:
$(a-1)(h-1)=(b-1)(e-1)$; $(c-1)(b-1)=(d-1)(g-1)$; $(e-1)(d-1)=(a-1)(f-1)$ and
$(g-1)(f-1)=(h-1)(c-1)$. This turns out to be the case.\\
If we are in case non-A, then the fourth point is $(1:1:0)$. In this case
we must have $x_1+x_3=0$, where we have the additional equations:
$b-a-h+e$, $d-c-b+g$, $f-e-d+a$ and $h-g-f+c$ are zero. Also this is can be verified.\\
\vskip 6pt
We will redo all the above
computations for the second case:
So now
$B=\{(a:b:1)$, $(c:d:1)$, $(e:f:1)$, $(g:h:1)\}$ and $C=\{(a:d:1)$, $(c:b:1)$, $(e:h:1)$, $(g:f:1)\}$. Equations
from the point $a_3$ (w.l.o.g.)
\[
 ah=be, ~~~ cf=dg, ~~~ ed=fa, ~~~ gb=hc.
\]

We start with case that $A$ contains $(1:1:1)$ as its fourth point:
In this case we have the additional equations: $(a-1)(f-1)=(b-1)(g-1)$, $(c-1)(h-1)=(d-1)(e-1)$,
$(e-1)(b-1)=(f-1)(c-1)$ and $(g-1)(d-1)=(h-1)(a-1)$.
We show that our matrix has rank $<7$, by giving an explicit relation between the columns.
\[
\begin{array}{l|rrrrrrr}
P     & x^2y & x^2z & y^2x & y^2z & z^2x & z^2y & xyz\\ \hline
(111) & 1    & 1    &  1  &  1   &  1   &  1   &  1\\
(ab1) & a^2b & a^2  & b^2a & b^2 &  a   &  b   & ab\\
(cd1) & c^2d & c^2  & d^2c & d^2 &  c   &  d   & cd\\
(ef1) & e^2f & e^2  & f^2e & f^2 &  e   &  f   & ef\\
(gh1) & g^2h & g^2  & h^2g & h^2 &  g   &  h   & gh\\
(ad1) & a^2d & a^2  & d^2a & d^2 &  a   &  d   & ad\\
(cb1) & c^2b & c^2  & b^2c & b^2 &  c   &  b   & cb\\
(eh1) & e^2h & e^2  & h^2e & h^2 &  e   &  h   & eh\\
(gf1) & g^2f & g^2  & f^2g & f^2 &  g   &  f   & gf
\end{array}
\]
We find that $\sum x_i\underline{c}_i=\underline{0}$ for the following coefficients:
\begin{eqnarray*}
x_1 & = & - b - d + f + h,\\
x_2 &=&    b d - f h,\\
x_3 &=&   a + c - e - g ,\\
x_4 &=& -  a c + e g,\\
x_5 &=&   b(f-d)(e+g) ,\\
x_6 &=&   e(c-g)(b+d) ,\\
x_7 &=& 0.
\end{eqnarray*}

What is left is the second case where the fourth point of $A$ is $(1:1:0)$.
Recall that $B=\{(a:b:1)$, $(c:d:1)$, $(e:f:1)$, $(g:h:1)\}$ and $C=\{(a:d:1)$, $(c:b:1)$, $(e:h:1)$, $(g:f:1)\}$, and
\[
   af=de, ~~~ ch=bg, ~~~ eb=ah, ~~~  gd=fc.
\]
From the point $(1:1:0)$ we get the additional equations:
$b-a=f-g$, $d-c=h-e$, $f-e=b-c$, $h-g=d-a$, and the whole system simplifies to:
$h=g+d-a$, $g=f-b+a$, $f=e+b-c$, plus the rest, let us first simplify these three:
$f=e+b-c$, $g=e-c+a$, $h=e-c+d$, and $af=de$, $ah=eb$, $gd=fc$.\\
So $f=e+b-c$, $g=e-c+a$, $h=e-c+d$, and $e(d-a)=a(b-c)$, $e(b-a)=a(d-c)$ and
$e(c-d)=-cd+da-cb+c^2$. Now the first leads to $(e+a)(b-d)=0$, so $e=-a$, and
by symmetry $g=-c$, $h=-b$ and $f=-d$ and again all is fine. What remains is to show that also these points
lie on a cubic.  We consider the matrix
\[
\begin{array}{l|rrrrrrr}
P     & x^2y-xy^2 & x^2z & y^2z & z^2x & z^2y & xyz\\ \hline
(110) & 0         & 0    &   0   &  0   &  0   &  0\\
(ab1) & a^2b-ab^2 & a^2  & b^2 &  a   &  b   & ab\\
(cd1) & c^2d-cd^2 & c^2  & d^2 &  c   &  d   & cd\\
(ef1) & e^2f-ef^2 & e^2  & f^2 &  e   &  f   & ef\\
(gh1) & g^2h-gh^2 & g^2  & h^2 &  g   &  h   & gh\\
(ad1) & a^2d-ad^2 & a^2  & d^2 &  a   &  d   & ad\\
(cb1) & c^2b-cb^2 & c^2  & b^2 &  c   &  b   & cb\\
(eh1) & e^2h-eh^2 & e^2  & h^2 &  e   &  h   & eh\\
(gf1) & g^2f-gf^2 & g^2  & f^2 &  g   &  f   & gf
\end{array}
\]
Now we use $e=-a$ etc.
The following matrix should have rank $<6$:
\[
\left(
\begin{array}{rrrrrrr}
 a^2b-ab^2 & a^2  & b^2 &  a   &  b   & ab\\
 c^2d-cd^2 & c^2  & d^2 &  c   &  d   & cd\\
 -a^2d+ad^2 & a^2  & d^2 &  -a   &  -d   & ad\\
 -c^2b+cb^2 & c^2  & b^2 &  -c   &  -b   & cb\\
 a^2d-ad^2 & a^2  & d^2 &  a   &  d   & ad\\
 c^2b-cb^2 & c^2  & b^2 &  c   &  b   & cb\\
 -a^2b+ab^2 & a^2  & b^2 &  -a   &  -b   & ab\\
 -c^2d+cd^2 & c^2  & d^2 &  -c   &  -d   & cd
\end{array}
\right)
\]
Since $e=-a$ the characteristic is not two.
We may simplify this matrix (after rearranging rows and
columns) to the direct sum of the two $4\times 3$ matrices
\[
\left(
\begin{array}{rrr}
a^2 &b^2 & ab\\
c^2 & d^2 & cd\\
a^2 & d^2 & ad\\
c^2 & b^2 & cb
\end{array}
\right)
\bigoplus
\left(
\begin{array}{rrr}
a^2b-ab^2 & a & b\\
c^2d-cd^2 &c &d\\
a^2d-ad^2 & a & d\\
c^2b-cb^2 & c & b
\end{array}
\right)
\]
and the only relation we have left is $d=a-b+c$. And now it can be verified
that the second matrix has rank $<3$.

\section*{Acknowledgments}
This research was done when the first
author was visiting the \it Seconda  Universit\`a degli Studi di Napoli \rm at
\it Caserta. \rm The authors wish to thank for their supports the research
group \it GNSAGA \rm of Italian \it Istituto Nazionale  di Alta Matematica \rm and
the \it Dipartimento di Matematica \rm of the \it Seconda Universit\`a degli Studi di Napoli. \rm

\vspace{0,5cm}\noindent {\em Authors' addresses}:

\vspace{0.2 cm} \noindent Aart BLOKHUIS \\
Department of Mathematics and Computing Science\\
Eindhoven University of Technology\\
P.O. Box 513
5600 MB Eindhoven\\
(The Netherlands).\\
E--mail: {\tt a.blokhuis@tue.nl}

\vspace{0.2cm}\noindent G\'abor KORCHM\'AROS\\ Dipartimento di
Matematica e Informatica\\ Universit\`a della Basilicata\\ Contrada Macchia
Romana\\ 85100 Potenza (Italy).\\E--mail: {\tt
gabor.korchmaros@unibas.it }

\vspace{0.2cm}\noindent Francesco MAZZOCCA\\ Dipartimento di Matematica\\
Seconda Universit\`a degli Studi di Napoli\\
Via Vivaldi 43 \\ 81100 Caserta (Italy)\\
E--mail: {\tt francesco.mazzocca@unina2.it }

\end{document}